\documentstyle[a4,11pt,amssymb]{article}

\def\section*#1{\begin{center} {\normalsize {\bf #1}} \end{center}}

\newcommand{\ee}{\hfill $\Box$}
\newcommand{\Proof}{\noindent {\sl Proof.} \hspace{3mm}}
\newcommand{\vv}{\vspace{2mm}}
\newcommand{\Frob}{\mbox{{\rm Frob}}}

\begin{document}

{\footnotesize \noindent Running head: Hyperelliptic class groups \hfill December 16, 1998

\noindent Math.\ Subj.\ Class.\ (1991): 11R29, 14H40 \hfill (slightly revised, August 25, 1999)}

\vspace{5cm}

\begin{center} 

{\Large The 2-primary class group of 

\vv

certain hyperelliptic curves}

\vv

{\sl by} Gunther Cornelissen

\vv

\end{center}

\vv

\vv

\begin{center} {\bf Introduction} \end{center}

\vv

\noindent  In a letter to Dirichlet, dated May 30, 1828 (\cite{Dirichlet}), Gau{\ss} considered the divisibility of the
class number of ${\bf Q}(\sqrt{-p})$ ($p$ prime, $=1$ mod 4) by 8.
Many variations on this theme can be found in the mathematical literature
of the subsequent centuries, either using quadratic forms  or class field theory ( R\'edei \cite{Redei:39}, Barrucand-Cohn \cite{Barrucand:69}, Hasse \cite{Hasse:69}, Kaplan \cite{Kaplan:73}, Stevenhagen \cite{Stevenhagen:93}). Of
these, the latter approach seems to give the most dense and structural arguments. 
To understand Gau{\ss}' letter from this point of view, one makes essential use of
the fact that the 2-primary part of the class group is cyclic, and that the 2-torsion in the
class group of ${\bf Q}(\sqrt{-p})$ is generated by the prime ideal above 2. In particular, the
norm of the ambiguous class is independent of $p$. The divisibility of the class number by 16 or even higher powers of 2 seems to be less tractable.

In the twenties, E.\ Artin in his thesis (\cite{Artin:24}) developed an 
arithmetic 
theory of quadratic extensions of the rational function field over a finite field, which shows a remarkable resemblance to the classical works of 
Gau{\ss} (quadratic forms), Dirichlet-Minkowski (units), Kummer-Dedekind (ideal class groups) and Riemann (zeta functions). With his terminology at hand, 
one can ask for the obvious analogue of the above: when is the class number
of (the ring of integers of) an imaginary hyperelliptic function field divisible by 2, 4 or 8? 

Throughout, we will assume that the characteristic of the ground field is different from 2.  The $2$-power rank of Artin-Schreier curves over fields of
characteristic 2 was studied by
van der Geer and van der Vlugt (\cite{Geer:91}). 

In view of arithmetic algebraic geometry, one has a lot of additional
geometric structure to study our question. The link between the class group 
and the Jacobian of the corresponding curve is well known (see the proof of the first corollary below). In particular, one can use the action of
Galois on the 2-power torsion in the Jacobian of the curve.  Such an approach
was undertaken in section 5 of \cite{Cornelissen:99} (following the suggestion of
a referee), and leads to a very satisfactory answer. But one can still wonder whether a class field theory approach to the problem is 
possible, and this is what we will attempt in this note. The geometric
result on the existence of rational torsion points on Jacobians of certain
hyperelliptic curves then comes out ``for free''. 
The main results can be stated as follows:

\vv

{\bf Theorem 1.} \ \ {\sl Fix a non-square $e \in {\bf F}_q$. For a prime $\frak p$ of  ${\bf F}_q[T]$, the class number of (the Dedekind ring) ${\bf F}_q[T,\sqrt{e{\frak p}}]$ is even if and only if $\deg {\frak p}$ is, and is divisible by 4 if and only if $\deg {\frak p}$ is. }

\vv

{\bf Theorem 2.} \ \ {\sl If $k$ is a fixed integer divisible by 4, then 
there exists a constant $C_{k}$ such that for all prime powers $q>C_{k}$ coprime to $k-4$, there exist two
primes ${\frak p}$ and ${\frak p}'$ of degree $k$ in ${\bf F}_q[T]$ and non-squares $e,e'$ in
${\bf F}_q$ such that the class numbers of ${\bf F}_q[T,\sqrt{e{\frak p}}]$ and ${\bf F}_q[T,\sqrt{e'{\frak p}'}]$ are different modulo 8.}

\vv 

{\bf Corollary 1.} \ \ {\sl Let $J(e{\frak p})$ be the Jacobian of a non-singular projective hyperelliptic curve whose affine equation is
$y^2=e{\frak p}$ for some non-square $e \in {\bf F}_q$  and
an irreducible polynomial $\frak p$ of degree $k$ in ${\bf F}_q[T]$. Then 
$J(e{\frak p})$ has an ${\bf F}_q$-rational two-torsion point if and only if $k$ is divisible by 4. Furthermore, for $4|k$ and all prime powers $q>C_{k}$ coprime to $k-4$ there are two
irreducible $\frak p$ and ${\frak p}'$ of degree $k$ and $e,e'$ non-squares
in ${\bf F}_q$ such that $J(e{\frak p})$ 
(resp.\ $J(e'{\frak p}')$) does (resp.\ does not) have an ${\bf F}_q$-rational 
point of exact order 4.}

\vv

{\bf Corollary 2.} \ \ {\sl Let $k$ be an integer divisible by 4.  For all prime powers $q>C_{k}$ coprime to $k-4$, there exists two
quaternion algebras over ${\bf F}_q(T)$ which are both ramified only at infinity
($T^{-1}$) and a unique finite prime of degree $k$ in ${\bf F}_q[T]$, but whose
type numbers (the number of non-conjugate maximal orders) have different parity.}

\vv

Contrary to the case of rational integers, the ambiguous class of
hyperelliptic function fields {\sl does} depend on the discriminant, and this 
turns out to be a severe obstruction to an immediate
translation of the classical argument. It also obscures the classical construction of a governing field for the 8-rank.

Recall that a Galois extension $\Omega_r$ of 
${\bf F}_q(T)$ is called a governing field for the $2^r$-rank with multiplier
$e$ if for all irreducible ${\frak p} \in {\bf F}_q[T]$, the class number of ${\bf F}_q(T,\sqrt{e{\frak p}})$ is divisible by $2^r$ if and only if $\frak p$ splits in $\Omega_r$ (with $\Omega_r$ independent of $\frak p$). From theorem 1, we see that $\Omega_2$ and $\Omega_4$  exist and are just the constant extensions ${\bf F}_{q^2}(T)$ and
${\bf F}_{q^4}(T)$. It is classically known that the 8-rank of 
${\bf Q}$-extensions with multiplier $-1$ is governed by ${\bf Q}(\sqrt[8]{-1},\sqrt{i+1})$; its construction depends heavily on the aforementioned independence ({\sl cf.\ }\cite{Stevenhagen:93} and the references therein). The analogue of this for ${\bf F}_q(T)$ is less clear: {\sl does there exist a field governing the 8-rank of class groups of hyperelliptic curves?} Let us only note that Bauer's theorem (``Galois extensions of number fields are determined by their splitting primes'', {\sl cf.} Neukirch \cite{Neukirch} VII.13.9) remains true for function fields up to constant extensions. This implies that $\Omega_r \otimes \bar{{\bf F}}_q \subseteq \Omega_{r+1} \otimes
\bar{{\bf F}}_q$.

The plan of this note is as follows: we recall the genus theory of Artin and
its connection to the parity of the class number. We then provide a
class field theory approach to 4-divisibility. However, we also give a second short argument using Drinfeld modular curves. Taken the construction of certain special ambiguous
classes for granted, we use class field theory to formulate
a criterium for 8-divisibility of the class number of curves corresponding to these classes. In the  next paragraph, discriminants having appropriate ambiguous classes are constructed; to produce them, we
construct certain ``lifts'' of the coefficients of $\frak p$ to the function field 
${\bf F}_q(t)$ and rely
on Chebotar\"ev's density theorem. 
Thus, the constant $C_k$ can be effectively estimated, and then the explicit construction
of $\frak p$ and ${\frak p}'$ is easy from the given data. The final paragraph
is devoted to the proof of the corollaries.

By similar constructions, it ought to be possible to 
surpress the divisibility conditions imposed on $q$ and $k-4$ 
in theorem 2. 

The results of this paper grew out of an attempt to get a better understanding of the interrelations between such class numbers, supersingular Drinfeld modules and Eisenstein series. For applications in that sense, see \cite{Cornelissen:99}.

\vv

\vv

\begin{center} {\bf 1. Notations -- Genus theory}  (E.\ Artin \cite{Artin:24}, section 11) \end{center}

\vv

\noindent Let ${\bf F}_q$ be a finite field with $q$ elements of characteristic $p \neq 2$, and let $K:={\bf F}_q(T)$ be the rational function field over ${\bf F}_q$
with maximal $T^{-1}$-order $A={\bf F}_q[T]$. Let $e$ be a non-square in 
${\bf F}_q$, and $\frak p$ an irreducible non-constant polynomial of degree $k$ in $A$. We write
$L=K(\sqrt{e
{\frak p}})$ for the quadratic extension of $K$ of discriminant $e {\frak p}$, and $\cal O$ for its ring of integers. 
$L$ is a so-called {\sl imaginary} quadratic extension of $K$ since $T^{-1}$ is
inert in $L$.
 
Let $h(e{\frak p})$ denote the class number of $\cal O$. Let $\cal C$ denote the 2-primary part of the divisor class group of the Dedekind ring
${\cal O}$.  Then $\cal C$ is a cyclic 2-group (since the discriminant has only one monic prime divisor), and it is trivial if and only if deg($\frak p$) is odd ({\sl cf.\ }Artin, section 11, Satz). This proves the first claim in theorem 1.

Artin also shows that if deg($\frak p$) is even, then the two-torsion of $\cal C$ is generated by
the ``ambiguous class'' $\cal A$, constructed as follows: consider the quadratic form
$$ \kappa (X,Y) = e X^2 + \alpha XY + Y^2, $$
where $\alpha \in {\bf F}_q$ is chosen such that $\kappa$ is irreducible over ${\bf F}_q$. If $(B,C) \in A$ are
such that 
$$ e{\frak p} =\kappa (B,C) \mbox{ with } \deg B < \deg C = \frac{k}{2}, $$
then in fact $\cal A$ is the class of the ideal $(C,B+\sqrt{e{\frak p}})$. A small computation shows that it has $K$-norm $N^L_K({\cal A}) = C$. As we have remarked in the introduction, this norm depends
on the discriminant of $L$. If $\frak p$ corresponds to $B$ and $C$ in this way, we will indicate this dependence by ${\frak p}(B,C)$, by slight abuse of notation. 

Let $H$ be the Hilbert class field of
$L$; this means the maximal abelian unramified extension of $L$ in which $T^{-1}$ is totally split (\cite{Hayes:79}). We will let ${\cal C}^i$ denote the groups $\{ c^i | c \in {\cal C} \}$. We will denote by $H_i$ the fixed field
of $H$ under ${\cal C}^i$. Because of the cyclic structure of ${\cal C}$, 
we see that $$2^i|h(e{\frak p}) \iff [H_{2^i}:L]=2^i.$$ 

\vv

\vv

\begin{center} {\bf 2. Proof of theorem 1} \end{center}

\vv

(2.1) {\sl Proof using class field theory.} \ \  The parity of the class number is given by the genus theory of \S 1. For the divisibility by 4, assume that $k=\deg({\frak p})$ is even, {\sl viz.}
$[H_{2}:L]=2$. The class number $h(e{\frak p})$ of $\cal O$ is divisible by 4 if and only
if $[H_4:H_2]=2$. 
We will use properties of the subfield $K(\sqrt{e})$, which is independent of 
${\frak p}$. This is actually the exact field of constants of $H_4$, so that in what follows we will not have to worry about unramified extensions of $K(\sqrt{e})$ in $H_4$. 
Since the ``genus field'' $H_2$
equals $L(\sqrt{e})$, $4$ divides $h(e{\frak p})$ if and only if $[H_4:K(\sqrt{e})] =4$. 

The extension $H_4/K(\sqrt{e})$ is unramified outside $\frak p$, and hence a subextension of the ray class field of $K(\sqrt{e})$ modulo $\frak p$. Let $R_{\frak p}$ be the ray class group of $K(\sqrt{e})$ at $\frak p$. Then Gal$(H_4/K(\sqrt{e}))$ is a quotient of $R_{\frak p}$. Since
$K(\sqrt{e})$ has class number one, it equals its own class field. If $\frak p$ splits as $\pi.\pi'$ in $K(\sqrt{e})$, then $R_{\frak p}$ is by class field theory equal to 
 $$ R_{\frak p} = [(A[\sqrt{e}] /{\pi})^* \times (A[\sqrt{e}] /{\pi'})^*]/{\bf F}^*_q(\sqrt{e}). \leqno{(2.1.1)} $$
(see Hayes, \cite{Hayes:79} \S 9). 
The situation is summarized in the following diagram:

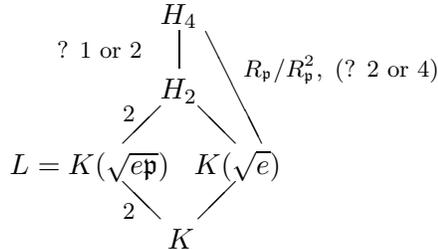
\begin{center}
\unitlength0.5mm

\begin{picture}(70,70)
\put(10,30){\line(1,1){10}}
\put(10,20){\line(1,-1){10}}
\put(30,40){\line(1,-1){10}}
\put(30,10){\line(1,1){10}}
\put(25,50){\line(0,1){10}}
\put(32,60){\line(1,-2){15}}
\put(-10,53){\mbox{ {\footnotesize ? 1 or 2}}}
\put(10,36){{\footnotesize 2}}
\put(42,48){{\footnotesize $R_{{\frak p}}/R_{{\frak p}}^2, \mbox{ (? 2 or 4)} $}}
\put(20,62){$H_4$}
\put(20,42){$H_2$}
\put(-20,22){$L=K(\sqrt{e{\frak p}})$}
\put(29,22){$K(\sqrt{e})$}
\put(22,2){$K$}
\put(10,10){{\footnotesize 2}}
\end{picture} \end{center}

\noindent  From this diagram, we see that the Galois group of the extension $H_4/K(\sqrt{e})$ is of exponent 2, and hence a surjective image of  $R_{\frak p}/
R_{\frak p}^2$.

If 4 divides $h(e{\frak p})$, then $[H_4:K(\sqrt{e})] =4$, so $R_{\frak p}/R^2_{\frak p}$ has order 4. By (2.1.1), this means that $\sqrt{e}$ is a square modulo 
$\pi$. But that happens if and only if $2|\deg(\pi)$, {\sl viz.,} $4|k$.

On the other hand, if $4|k$, then $R_{\frak p}/R^2_{\frak p}$ has order 4. To
it corresponds an extension of $K(\sqrt{e})$ with Galois group ${\bf Z}/2 
\times {\bf Z}/2$ by (2.1.1), which is only ramified at the primes above $\frak p$, at
most with ramification index 2. Hence it is contained in $H_4$. But then $[H_4:K(\sqrt{e})] \geq 4$, {\sl i.e.,} 4 divides $h(e{\frak p})$.  \ee

\vv

(2.2) {\sl Proof using modular curves.} \ \ Let $k$ again be even. In \cite{Gekeler:86}, 
E.-U.\ Gekeler shows that 
$$ 4 g({\frak p})  =  2 \frac{q^k-1}{q^2-1} + h(e{\frak p}), $$
where $g({\frak p})$ is an integer. Let it suffice for the {\sl cognoscenti} to 
remark that $g({\frak p})$ is the genus of the quotient of the Drinfeld modular curve 
$X_0({\frak p})$ by the $\frak p$-Atkin-Lehner involution, or equivalently, the
number of quadratic $j$-invariants for rank-two Drinfeld modules that are
supersingular modulo $\frak p$.
 
Since it is elementary to see that $2 (q^k-1)/(q^2-1) =  k \mbox{ mod } 4$, we again find that $4|h(e{\frak p})$ if and only if $4|k$. \ee

\vv

\vv

\begin{center} {\bf 3. Proof of theorem 2}  \end{center}

\vv

(3.1) The cyclic structure of $\cal C$ implies that $h(e{\frak p})$ is 
divisible by 8 if and only if ${\cal A} \in {\cal C}^4$. We will now proceed
to construct, for $k$ divisible by 4 and $q$ large enough, a pair $(B,C)$ 
of a particular form, and derive a criterium for the corresponding field $K(\sqrt{e {\frak p}(B,C)})$ to have class groups of prescribed 8-rank. {\sl Throughout this section, we set $l=k/2$} (which is even), and use the notations of the first
paragraph.  It will also be useful to keep in mind the diagram
of \S 2. 

\vv

(3.2) Let $C=T^{l-2} Q(T)$ for some quadratic irreducible polynomial
$Q(T)=T^2+aT+b$ over ${\bf F}_q$, and let $B= b_0 T^{l-2}+ b_1 \in A$. Write $L=K(\sqrt{e{\frak p}(B,C)})$. It is easy to see that all factors of
$C$ split in $L/K$, since $e{\frak p} = B^2 \mbox{ mod } C$ and the conductor
of $A[e{\frak p}]$ is trivial. Hence we can write $T={\cal T} {\cal T}'$ and $Q={\cal Q} {\cal Q}'$ in $L$.  Then the ambiguous class  is ${\cal A} = {\cal T}^{l-2} {\cal Q}$. 

Since $T$ is of degree one, it is not split in the (constant) extension $K(\sqrt{e})={\bf F}_{q^2} \otimes K$.  Since the ramification index and residue class degree of $T$ in $H_2/K$ are the same whether computed via $L$ or via
$K(\sqrt{e})$, we find that $\cal T$ is not split
in $H_2/L$ either, {\sl i.e.}, $\Frob_{\cal T} \neq 1$ in Gal$(H_2/L) = {\cal C}/{\cal C}^2$. On the other hand, $Q$ splits in $K(\sqrt{e})/K$, say as
$Q={\cal L} \cdot {\cal L}'$, and one sees that $\cal L$ and ${\cal L}'$ also split in $H_2=L(\sqrt{e})/K(\sqrt{e})$ (again since $e{\frak p}$ is a square
modulo ${\cal L}$). Hence by a similar argument as before, ${\cal Q}$ splits in $H_2/L$, {\sl i.e.,} $\Frob_{\cal Q} = 1$ in Gal$(H_2/L)$. 
If we let $\sigma$ denote a generator of the cyclic group Gal$(\lim\limits_{\longrightarrow} H_{2^i}/L)\cong {\cal C}$, then
we can write $\Frob_{\cal T}^2=\sigma^{4m+2}$, and $\Frob_{\cal Q}=\sigma^{2n}$ for some $m,n \in {\bf Z}$. 
We now distinguish two cases:

\vv

(3.2.1) {\sl First case: $l-2$ is divisible by 4.} \ \ The class number $h(e{\frak p})$ is divisible by $8$ if and only if ${\cal A} \in {\cal C}^4$, {\sl i.e.}, $\Frob_{\cal T}^{l-2} \circ \Frob_{\cal Q}=1$ in Gal$(H_4/L)$. Since the latter group is of exponent 4, this is equivalent to Frob$_{\cal Q}$ acting trivial in Gal$(H_4/L)$. So we want ${\cal Q}$ to split completely in $H_4/L$. But this is equivalent
to $\cal L$ splitting completely in $H_4/K(\sqrt{e})$ (then the same follows for ${\cal L}'$). Using the description of its Galois group in terms of the ray
class group at $\frak p$ given earlier, we want that
Frob$_{\cal L}$ acts trivial in $R_{\frak p}/R^2_{\frak p}$. Writing ${\frak p}=\pi \pi'$ in $K(\sqrt{e})$
as before, and using (2.1.1),  this
is equivalent to $\cal L$ being a square modulo $\pi$ in ${\bf F}_{q^2}$ (the same then immediately holds modulo $\pi'$ since ${\frak p}$ is a square modulo
${\cal L}$). 

We can reformulate this criterium using quadratic residue symbols $(\frac{\cdot}{\cdot})$ for ${\bf F}_{q^2}(T)$ as follows. Let us factor the quadratic form $\kappa$ 
over ${\bf F}_{q^2}$ as 
$$ \kappa (X,Y) = e(X-\delta Y)(X-\bar{\delta} Y) \mbox{ over } {\bf F}_{q^2}, $$
where $2 e \delta = \alpha  \pm \sqrt{\alpha^2-4e}$.
Then $\pi=C-\delta B= T^{l-2} Q(T) - \delta B(T)$ for some choice
of $\delta$. Let $\lambda$ be an root of $Q$ in ${\bf F}_{q^2}$, say, ${\cal L}=T-\lambda$. Our criterium reads
\begin{eqnarray*} 8 | h(e{\frak p}) &\iff& (\frac{{\cal L}}{\pi})=1 \iff (\frac{\pi}{{\cal L}})=1\\ &\iff&
\delta B(\lambda) = \mbox{ square in }  {\bf F}_{q^2}, \end{eqnarray*}
using quadratic reciprocity and the fact that $\cal L$ is a factor of $Q$. 

\vv

(3.2.2) {\sl Second case: $l$ is divisible by 4.} \ \ Now,  the class number $h(e{\frak p})$ is divisible by $8$ if and only if $\Frob_{\cal T}^2 \circ \Frob_{\cal Q}=1$ in Gal$(H_4/L)$, {\sl viz.,} $\sigma^{4m+2+2n} = 1$ in Gal$(H_4/L)$. Since the latter group is of exponent 4, this happens
if and only if $n$ is odd, {\sl i.e.}, $\Frob_{\cal Q}\neq 1$ in Gal$(H_4/L)$. We can then follow
the argument of (3.2.1) to see that this is equivalent to
$$ 8 | h(e{\frak p}) \iff \delta B(\lambda) \neq \mbox{ square in }  {\bf F}_{q^2}. $$

\vv

(3.3) For the above constructions to work, we have to require additionally
that ${\frak p}(B,C)$ is irreducible over ${\bf F}_q$, and that $e$ is a non-square in ${\bf F}_q$. Factor $\kappa$
as before. The fact that $\kappa(0,1)=1$ implies that $\delta \bar{\delta} = e^{-1}$. In the next section we
will prove the following proposition: 

\vv

{\bf (3.4) Proposition.} \ \ {\sl Fix an even integer $l$. There is a constant 
$C_l$ such that for all prime powers $q>C_l$ coprime to $l-2$, there
exists a quadratic irreducible $Q$ over ${\bf F}_q$, $B^\pm(T)= b^\pm_0 T^{l-2} + b^\pm_1 \in {\bf F}_q[T]$ and  $\delta^\pm 
\in {\bf F}_{q^2}$ such that 
$$ T^{l-2}Q(T)-\delta^\pm B^\pm(T)$$ is irreducible over ${\bf F}_{q^2}$, and 
such that  
$B^- (\lambda)$ is (and $\delta^{\pm}$ and $B^+(\lambda)$ are not) a square in  ${\bf F}_{q^2}, $
where $\lambda$ is a root of $Q$ in ${\bf F}_{q^2}$.}

\vv

\noindent Let us show how this leads to our requirements. Remark that once $\delta^\pm$ is given, it is easy to compute a corresponding
$\alpha$. Since $\delta^\pm$ is not a square in ${\bf F}_{q^2}$, the inverse of its norm, $e$, is not a square in ${\bf F}_q$. Also, since $\delta^\pm \notin {\bf F}_q$, the corresponding $\frak p$ is irreducible. Namely, if $C-\delta^\pm B^\pm \in {\bf F}_q[T]$, then it would factor over ${\bf F}_{q^2}$ since
its degree $l$ is even.  \ee

\vv

\vv

\begin{center} {\bf 4. Construction of appropriate discriminants} \end{center}

\vv

{\bf (4.1) Proposition.} \ \ {\sl Let $q$ be an odd prime power and $l$ an
even integer such that $l-2$ is coprime to $q$. Let $Q(T)=T^2+aT+b$ be an irreducible quadratic polynomial over ${\bf F}_q$ with $a \neq 0$. Let $B(T)= b_0 T^{l-2} + b_1$ be an irreducible polynomial in ${\bf F}_q[T]$ with $b_0,b_1 \neq 0$. Assume that $B(\lambda) \neq 0$ for the roots $\lambda$ of $Q$. Then the polynomial $$f:=T^{l-2} Q(T) - t B(T) $$
has Galois group $S_l$ over ${\bf F}_{q^2}(t)$ and  its splitting field has
${\bf F}_{q^2}$ as
its exact field of constants.}

\vv

\Proof We see that $f$ is irreducible over ${\bf F}_p(t)$ using Gau{\ss}' lemma since it is irreducible as
a (linear) polynomial in $t$. Let $G$ be its Galois
group over ${\bf F}_{q^2}(t)$.

We will first prove that $G$ is primitive. Remark that it suffices to show that $G$ is 2-transitive (\cite{Wielandt:Werke}, Theorem 9.6). For this, it suffices to show that the
stabilizer of any root $\alpha$ of $f$ is transitive. We now appeal to the following ``twisted derivative''-trick of Abhyankar's (\cite{Abhyankar:92}, \S  18), which says that this stabilizer
is 
$$ \mbox{Gal}(F/{\bf F}_{q^2}(\alpha)), \mbox{ where } F=\frac{f(T) - f(\alpha)}{T-\alpha}. $$
If we compute $F$ using $t=\alpha^{l-2} Q(\alpha) B(\alpha)^{-1}$, we get 
$$ F=T^{l-1} + (\alpha+a) T^{l-2} + b_1 \frac{Q(\alpha)}{B(\alpha)} \sum_{i=0}^{l-3} \alpha^i T^{l-3-i}. $$
All we have to show is that $F$ is irreducible over ${\bf F}_{q^2}(\alpha)$, which is a  rational function field 
(remark that $\alpha$ is transcendental over ${\bf F}_{q^2}$). The Newton polygon for the valuation corresponding to the prime factors of $Q(\alpha)$ in ${\bf F}_{q^2}[\alpha]$ contains a straight line segment from
$(0,1)$ to $(l-2,0)$ which goes through no integer lattice points. Hence
if $F$ is reducible, it has a root, say $T_0$, which is not divisible by any factor
of $Q(\alpha)$. By assumption, $B(\alpha)$ has at most two factors over ${\bf F}_{q^2}[\alpha]$, and the Newton polygon of $F$ for such a factor is a straight line from $(0,-1)$ to $(l-2,1)$ followed by a segment of slope one. Hence such a factor of $B(\alpha)$ occurs at most once in the factorization of $T_0$. From the Newton polygons of $F$ for all other
finite valuations, one sees that the only possible further divisor of $T_0$ is $\alpha$, with valuation $-1$. 
The Newton polygon of $F$ for the valuation $-\deg_\alpha$ is a straight line
from $(0,-1)$ to $(l-2,-1)$ having slope zero, followed by a segment from $(l-2,-1)$ to $(l-1,0)$
having slope one. Since $F$ has no constant roots, we find that $\deg_\alpha
T_0 = 1$. Suppose first that $l \neq 4$. Then a factor of $B(\alpha)$ (which has degree $\frac{1}{2}(l-2)>1$) cannot divide $T_0$. The data just computed imply that $T_0$ is constant, a contradiction. Similarly, if $l=4$, let $\beta_1,\beta_2$ denote the two (linear) factors of $B$ over ${\bf F}_{q^2}$. Then the only possibilities are that $T_0$ is a scalar multiple of $\beta_1(\alpha), \beta_2(\alpha)$ or $B(\alpha)\alpha^{-1}$. If we expand the two leading terms of $F(T_0)$ in $\alpha$ in each of these cases, we find that $\beta_1,\beta_2$ are defined over ${\bf F}_q$, which contradicts the supposed irreducibility of $B$ over ${\bf F}_q$. 
\vv

We will now prove that $G$ contains a $2$-cycle. Make the following change
of variables: $Y=1/T,u=1/t$. Then 
$$ f = b_1 Y^l + b_0 Y^2 - u (b Y^2 + a Y + 1). $$
This is an Eisenstein polynomial for the prime $u$, and its reduction is
$f = Y^2 ( b_1 Y^{l-2} + b_0 ) \mbox{ mod } u. $ 
Since $ b_0 Y^{l-2} + b_1$ has no multiple roots over ${\bf F}_{q^2}$, the factorization of $f$ over the $u$-completion of ${\bf F}_{q^2}(t)$ consists
of an Eisenstein polynomial of degree $2$ multiplied by a polynomial without multiple factors mod $u$. Proposition (3.1) in \cite{Cornelissen:99} says that the inertia group of the splitting field of an Eisenstein polynomial over
a local field whose degree, say $N$, is prime to the residue characteristic contains an $N$-cycle. If we apply it to our situation, we find that the inertia group
of $1/t$ in $G$ contains a 2-cycle.

We now appeal to a result of Jordan (\cite{Neumann:85}) which says that a
primitive permutation group of degree $l$ containing a $2$-cycle is
$(l-1)$-transitive. Hence $G=S_l$. 

If ${\bf F}_{q^{2N}}$ is the
exact field of constants of the splitting field ${\cal F}$ of $f$ over ${\bf F}_{q^2}$, then ${\bf F}_{q^{2N}}(t)$ is the fixed
field of the group $G'$ generated by all inertia groups. Since $A_l$ is simple
for $l>4$, we find that $G'$ is either $A_l$ or $S_l$. But above we
have constructed an even element in such an inertia group, so this proves that $N=1$ in those cases. If $l=4$, then the only other normal subgroup of $A_l$
is generated by products of two transpositions, but this group can also
not equal $G'$ since the latter contains a transposition.   \ee

\vv

(4.2) {\sl Proof of proposition (3.4).} \ \ 
 Choose $Q=T^2+aT+b$ irreducible over ${\bf F}_q$ with $a \neq 0$, having a root $\lambda$ which is a multiplicative generator for 
${\bf F}^*_{q^2}$. Choose $q$ big enough to have $q+1>l-2$. Then $\lambda^{l-2} \notin {\bf F}_q$. This implies that the set $\{ b_0 \lambda^{l-2} + b_1 : b_0,b_1 \in {\bf F}_q \}$ is the whole of ${\bf F}_{q^2}$. Hence we can choose
 $B^\pm$ such that $B^-(\lambda)$ is (and $B^+(\lambda)$ is not) a square
in ${\bf F}_{q^2}$ (with the corresponding $b_i^\pm \neq 0$) by imposing one linear relation between $b_0$ and $b_1$. Also, one can choose such $b^\pm_0,b^\pm_1$ such that $B^{\pm}(T)$ is irreducible over ${\bf F}_q[T]$
(this is the easy case of the Hansen-Mullen conjecture, {\sl cf.\ }\cite{Hansen:92}). Let $f^\pm = T^{l-2}Q(T)-tB^\pm(T)$. By the previous 
proposition, we find that the splitting field ${\cal F}^\pm$ of $f^\pm$ over
${\bf F}_{q^2}(t)$ has an $l$-cycle $\sigma_l$ in its Galois group, and ${\bf F}_{q^2}$
as its exact field of constants. 

There are now two possibilities: either the Galois group of $f^\pm$ over
${\bf F}_{q^2}(\sqrt{t})$ is $A_l$ or $S_l$, depending on whether the discriminant
of $f^\pm$ equals $t$ up to squares or not. Suppose it is $A_l$. 
We can apply Chebotar\"ev's theorem (as in the appendix of Geyer and Jarden \cite{Geyer:98}) to  the Galois extension ${\cal F}^\pm$
of ${\bf F}_{q^2}(t)$ to find that for big enough $q$, there exists
primes $P^\pm=t-\delta^\pm$ of degree one in  ${\bf F}_{q^2}(t)$ whose Frobenius elements act like the $l$-cycle $\sigma_l$ on the roots of $f^\pm$.
Since ${\bf F}_{q^2}(\sqrt{t})$ is the fixed field of $A_l$ in ${\cal F}^\pm$,
and the $l$-cycle $\sigma_l$ is even, it acts non-trivially on $\sqrt{t}$, and hence
the same holds for the Frobenius elements of $P^\pm$.  This means that $t=\delta^\pm$ is never a square modulo $P^\pm$ in ${\bf F}_{q^2}$. On the other
hand, $B^+(\lambda)$ is, and $B^-(\lambda)$ is not, a square in ${\bf F}_{q^2}$.

If the Galois group of $f^\pm$ over ${\bf F}_{q^2}(\sqrt{t})$ is $S_l$, then
the fields ${\cal F}^\pm$ and ${\bf F}_{q^2}(\sqrt{t})$ are disjoint over 
${\bf F}_{q^2}(t)$, and hence the Galois group of ${\cal F}^\pm(\sqrt{t})$ over ${\bf F}_{q^2}(t)$ equals $S_l \times {\bf Z}/2$, where the generator of ${\bf Z}/2$ acts like $\sqrt{t} \rightarrow -\sqrt{t}$. 
We can then apply Chebotar\"ev's theorem to the extension ${\cal F}(\sqrt{t})$ of ${\bf F}_{q^2}(t)$ to find primes $P^\pm = t-\delta^\pm$ of degree one in ${\bf F}_{q^2}(t)$ whose Frobenius elements act like $\sigma_l \times (-1)$ in $S_l \times {\bf Z}/2$. This means again that $\delta^\pm$ are non-squares in ${\bf F}_{q^2}$ and
the corresponding statements about $B^\pm(\lambda)$ are satisfied. 
 This proves
the proposition. \ee

\vv

\vv

\begin{center} {\bf 5. Proof of the corollaries} \end{center}

\vv

(5.1) {\sl Proof of corollary 1.} \ \ From genus theory, we know that 
${\cal C}$ is cyclic. The result then follows immediately from the well-known exact sequence 
$$ 0 \rightarrow J[2^\infty]({\bf F}_q) \rightarrow {\cal C} \rightarrow {\bf Z}/ d {\bf Z} \rightarrow 0, $$
where $d \in \{1,2\}$ is the degree of $T^{-1}$ in $K(\sqrt{e{\frak p}})$ (so 
$d=2$ if and only if $k$ is even). \ee

\vv

(5.2) {\sl Proof of corollary 2.} \ \ This follows from the formula in (2.2) since the type
number $t({\frak p})$ of the quaternion algebra ramified at the finite prime $\frak p$ of degree $k$ equals $g({\frak p})$ (Gekeler, {\sl loc.\ cit.}). Actually, for $4 | \deg {\frak p}$,
$$ 2 | t({\frak p}) \iff \left \{ \begin{array}{ll} 8 | h(e{\frak p})  & \mbox{ if } 8 | k \\ 8 \hspace{-2mm} \not | h(e{\frak p})  & \mbox{ if } 8 \hspace{-2mm} \not | k \end{array} \right. , $$
as follows from the formula for $g({\frak p})$ given in (2.2). \ee

\vv

\vv

{\footnotesize

\noindent {\bf Acknowledgments.} The author is post-doctoral fellow of the 
Fund for Scientific Research - Flanders (FWO -Vlaanderen). This work 
was done while visiting the MPIM. 

\vv

\vv

\bibliographystyle{amsplain}
\providecommand{\bysame}{\leavevmode\hbox to3em{\hrulefill}\thinspace}

\vv

\vv

\noindent Max-Planck-Institut f{\"u}r Mathematik, Gottfried-Claren-Stra{\ss}e 26,
D-53225 Bonn {\sl (current)}

\vv

\noindent University of Gent, Dept.\ of Pure Mathematics, 
Galglaan 2, B-9000 Gent 

\vv

\noindent e-mail: {\tt gc@cage.rug.ac.be}
}

\end{document}